\input amstex

\documentstyle{amsppt}
\overfullrule=0pt
\magnification = \magstep 1
\hsize=6.5 true in
\vsize=8.9 true in

\topmatter

\title 
On the Rationality of $SU(r,d)$ 
\endtitle

\author David C. Butler \endauthor


\address School of Mathematics, 
Tata Institute of Fundamental Research, 
Homi Bhabha road, Bombay, India  \endaddress

\address e-mail: butler\@math.tifr.res.in \endaddress

\endtopmatter

\document

\baselineskip = 20 pt
\lineskip = 11 pt
\lineskiplimit = 10 pt

\define\aut{\operatorname{Aut} }
\define\oc{\Cal O_C}
\define\ext{\operatorname{Ext}^1 }
\define\prj{ {\Bbb P} }

\define\cplus{\overset c \to \oplus}

\define\rank{\operatorname{rank}}
\define\mup{\mu^{\max } }
\define\mudn{\mu^{\min} }
\define\mdn{\mu^{\min } }
\define\bho{H^0(C,}
\define\ho{h^0(C,}
\define\bhi{H^1(C,}
\define\hi{h^1(C,}
\define\gr{\Bbb G}
\define\alph{\frac{\alpha}{2} }


\heading Introduction \endheading

Let $C$ be a curve of genus $g$ over the complex numbers.
For $0< d< r$, let $SU=SU(r, r(g-1)+d)$ be the moduli space of bundles 
with rank $r$ and fixed determinant of $\deg r(g-1) + d$ over $C$.
This paper shows that if $d$ divides $(r \pm 1)$, 
then $SU$ is rational (no matter what the value of $g$).
In fact if the prime divisors of $\delta >0$  are divisors
of $r$, and $d$ divides $r-\delta>0$ then $SU$ is rational
for any genus.  We conjecture that this holds for
$r+\delta$ if $0 < \delta < r$.

Research on the  rationality of $SU$ began 
shortly after the construction of the Moduli space. 
In 1964, Tyurin \cite {9}
proved the case of rank 2 bundles of odd degree.  He went on
to claim rationality for all ranks and degrees \cite {10},
but the paper contained fatal errors.  Newstead corrected
some errors, and in his unpublished 
thesis, he proved the case $d=1$ for any $r$. 
He then went on to claim a proof for $(r,d)=1$ \cite{6}.
That proof contained a trivial but fatal error.
He concluded in his correction \cite {6} that  for any pair
$(r,d)$ which are coprime the theorem
holds for infinitely many genera.  Which suggests, but does
not prove, that the theorem is true.  Ballico then noted that 
Newstead's argument  proved $SU$ was stably rational \cite {1}, 
where X is stably rational with
level $k$ if for all $n \ge k$, $X \times \prj^n$ is rational.
Unfortunately, Beauville, Colliot-Th\'el\`ene,
Sansuc and Swinnerton-Dyer proved stable rationality
does not imply rationality \cite {2}.  
However, Ballico noted that 
stable rationality with small level and Newstead's arguments
would  yield strong results.  Aware of this,
Boden and Yokogawa \cite{3}, used parabolic bundles to
prove $SU(r,d)$ is stably rational of level $\le r-1$.
This, and the argument of Newstead proved that if
$(r,d)=1$, then $SU$ is rational provided either
$(g,d)=1$ or $(g,r-d)=1$.  And this result generalizes 
Newstead's results.

Now our results continue the theme of widening the
number of genera for which rationality is known.
In particular, we find infinitely many pairs
$(r,d)$ with $(r,d)=1$ where $SU(r,d)$ is rational.

\proclaim {Main Theorem} 
Let C be a smooth projective irreducible curve 
over $\Bbb C$.
Then $SU(r,d^{'})$ is rational if $d^{'}\equiv \pm d \pmod r$, $(r,d)=1$ and if
\roster
\item "A)" $d$ divides $(r \pm 1)$, or if
\item "B)" $d$ divides $(r-\delta)$, where 
for all primes $p$ such that $p|\delta$, $p|r$, and $0<\delta<r$.
\endroster
\endproclaim

Along the way we wind up proving a second, seemingly
unrelated result:

\proclaim {Second Theorem} 
Let C be as above, and let $\pi: \tilde C \to C$ 
be a cyclic etale cover of deg $e$.  
For a general stable bundle $E$ with $\mu (\pi^* E) > g(\tilde C)$, 
$\pi^*(E)$ is generated by global sections.
\endproclaim

Before we outline the proof, it should be noted that
$SU(r,d) \cong SU(r,rk\pm d)$ for any integer $k$
by tensoring and/or dualizing.

Now we outline the proof giving a historic perspective.
It begins with Newstead who observed that since a general
bundle with rank $r$ and $\deg r(g-1) + d$ had $h^0 =d$,
there is a unique subbundle $I=\overset d \to \oplus \Cal O_C$
to a general $E$.
This gives an exact sequence
$$0\to I \to E \to F\to 0.$$
It seems likely, and turns out to be true that for general $E$,
$F$ is a stable bundle, and so a general $E$ is given by an extension of a
$\rank r-d=n$, $\deg r(g-1) +d$ vector bundle $F$ by $I$.
Now we know $SU=SU(r,r(g-1)+d)$ has a map to $SU^{'}=SU(r-d, r(g-1)+d)$
whose fiber over a closed point $[F]$ is 
$\ext (F,I)\cong \bho \omega_C\otimes F)^{*}$ modulo
the action of $\aut (I)$, and this turns out to be a Grassmannian.
But even if $SU^{'}=SU(r-d, r(g-1) +d)$ is rational, $SU$ need not be rational.
Because unfortunately, if $X$ maps to a rational variety with rational
fibers it is not necessarily true that $X$ is rational.
But if $(r-d, r(g-1)+d)=1$, then  $SU^{'}$ has a Poincare bundle
$\Cal P$ on $SU^{'} \times C$.  Tensoring by the pullback of $\omega_C$
then pushing down to $SU^{'}$ and taking $d$ direct summands, we get a Zariski
open set with a bundle parameterizing $\ext (F,I)$.  Taking $\gr$ the 
quotient by $\operatorname{Aut}(I)$, we get a parameterization of an open 
set of $SU$ by a Grassmannian bundle over $SU^{'}$.  
Since it is a bundle, it is trivial 
(on a perhaps smaller set) and so $SU$ is birational to $SU^{'}\times \gr$.
If $SU^{'}$ is rational by inductive hypothesis
then $SU$ is rational.  Unfortunately there is one problem.
$r(g-1) +d$ and $r-d$ need not be relatively prime. And if they
are not, there is no Poincare bundle.

To see that the above happens consider the following example.
Set $r=5$, $g=6$ and $d=2$.  $5(6-1)+2=27$ and is not relatively
prime to $(5-2)=3$.  The same thing happens for $d=5-2=3$.
Now consider another example.  Set $r=7$, $g=6$ and $d=2$.
$7(6-1)+2=37$ which is relatively prime to $(7-2)$.
But $SU^{'}=SU(5, 5(6-1)+2)$ which is the previous example.
We can try $d=7-2=5$ but $7(6-1)+5=40$ is not relatively
prime to $7-5=2$.  The whole inductive argument breaks
down unless $g$ is chosen carefully, or if $d=\pm 1$, 
in which case substituting
$d=r-1$ (By tensoring and perhaps dualizing) gives a
rank 1 cokernel which is rational because it is a point 
and there is only one fiber.

Now Ballico notes that in the coprime-prime case 
we do not need the $SU^{'}$ to be rational, but for some
product of a rational variety with $SU^{'}$ to be rational.
If $SU^{'}$ is stably rational with level less than or equal
to the dimension of the Grassmannian, then $SU$ is rational.
Unfortunately Ballico could not prove the moduli spaces
were stably rational of small enough level; although
he did prove them stably rational.

Boden and Yokogawa then proved rationality for most
parabolic bundles.  As a corollary they proved
that bundles with $(\rank, \deg)=(r,d)=1$ are stably rational
with level $r-1$.  This implies $SU$ is stable if
$(r(g-1),r-d)=1$ or $(r(g-1),d)=1$.  And this simplifies
to $(g,d)=1$ or $(g,r-d)=1$.

Now the author enters the picture.  It would seem at first that
Newstead's original argument had reached its logical
conclusion.  But it could be improved by changing its
very beginning.  Newstead begins with the trivial bundle $I$
of $\rank d$, because its embedding into a general $E$
with $[E]\in SU$ is unique up to automorphisms of $I$.
But there are other bundles with this property.  Namely
any $\deg 0$ twist of $I$ or any direct sum of $\deg 0$ twists of $I$.
We use the latter.  Let $\{L_1, \dots L_e\}$ be a set
of mutually non-isomorphic $\deg 0$ line bundles and let $U_i$ be a rank
$d$ vector space.  Then set $V_i=U_i\otimes L_i$ and
$W=\overset {e} \to \oplus V_i$.
The bulk of our proof is to show that for any $W$ 
with $ed < r$ and a general $E$ 
the sequence $$0\to W \to E \to F\to 0$$
has a stable bundle $F$ as  cokernel. 
And when $\rank W > \rank E$ we have a surjection
of $W$ to $E$ giving an exact sequence,
$$0\to K_{E,W} \to W \to E \to 0.$$
For the plus case in part A, $K_{E,W}$ is a fixed line bundle
we call $A$.  And $SU$
is isomorphic to  $\bho  A^* \otimes W)/\aut (W)$,
and this turns out to be a product of Grassmannians.

Now for part B (which implies the minus case in part A).  
If $ed= r-\delta$
where every prime factor of $\delta$ is a prime factor of
$r$, then $(\delta, r(g-1) +d)=1$ since no divisor of
$r$ divides $d$.  That means $SU^{'}=SU(\delta, r(g-1)+d)$
has a Poincare bundle.
A general $E$ is an extension of $F\in SU^{'}$ by $W$
And since $SU^{'}$ has a Poincare bundle, it has a bundle
whose  fibers are $\ext(F_t,W) \cong \overset {i} \to \oplus
\bho \omega_C \otimes L_i^* \otimes F_t)^* \otimes U_i$
Taking an affine subset of $SU^{'}$ we can form a universal
extension space.  This gives us a map from the bundle
over $SU^{'}$ to $SU$ which is dominant and has fibers
$\aut (W)$.  We take the quotient by the automorphism
group and get a birational map from $SU^{'}$ cross a product
of Grassmannians to $SU$.  Since the product of Grassmanians
has dimension $\ge r-1$, $SU$ is rational.

It would be nice to have similar results for when
$ed=r+\delta$ but they are held up because we
do not know when $K_{E,W}$ is stable.  In the
final section we make a natural conjecture
and offer some evidence.  We will also prove
the Second Theorem.

I thank Johns Hopkins University and the Tata Institute for
Fundamental Research for their hospitality.  I would like
to thank Shokurov who shared his expertise in birational
geometry.  I would also like to thank Mehta, Nitsure, Ramadas
Ramanan and Srinivas for answering questions so numerous that I do not
know who answered what.  But I do know they were answered.

\heading \S 1 Preliminaries \endheading

$C$ is a smooth projective irreducible curve of genus $g$
over $\Bbb C$

For a vector bundle $E$, $\mu(E)= \deg(E)/\operatorname{rank}(E)$. 

$E$ is stable (semistable) if $\mu(S) < \mu(E)$ $(\le)$
for all proper subbundles. 
If $E$ is not semistable, then it is called unstable.
If $E$ is not stable, then it is called nonstable.
If $E$ is semistable but not stable, then it is called strictly semistable.

As in the introduction $SU = SU(r, r(g-1)+d)$ and $U= U(r, r(g-1) +d)$,
the moduli spaces of stable bundles of $\operatorname{rank} r$
and $\deg r(g-1)+d$ with fixed and unfixed determinants
respectively.  It is assumed that $0< d < r$.

By abuse of notation we say $E\in SU$ (or $U$)
if the isomorphism class of the vector bundle is in
the moduli space.

We will use $E$ to denote a fixed bundle and $E_t$
to denote a general bundle in the moduli space.

$L$ will denote a single linebundle of $\deg =0$.

We will consider sets of mutually non-isomorphic line bundles of $\deg 0$:
$\{L_1, L_2, \dots L_e\}$.  Our results will not
depend on the the choice (except for the Second Theorem).

$U_i = \overset d \to \oplus k$ where $k$ is the ground field.

$V_i = U_i\otimes L_i$.

$W_k = V_1 \oplus \dots \oplus  V_k$.

$W=W_e$.
 
We will use dimension counts to prove our results.
These arguments are hampered because many bundles do
not live on a moduli space, and even those that do
may not have a universal or Poincare bundle.  But we
may (after Sundaram \cite {8}) form very useful parameter
spaces of bundles.

\subheading{Definition}
Let Cl be a class of vector bundles over $C$.
A parameter space $T$ of Cl is a variety
along with a bundle $\Cal E$ over $T\times C$
such that any $E$ in Cl is isomorphic to
the restriction $E_t$ of $\Cal E$ to some ${t}\times C$.

\subheading {Remark 1}
It is not assumed that  $E_t \cong E_s$ implies
$s=t$.  And even if it does, 
parameter spaces do not generally
have the universal properties of a moduli space.

Stable vector bundles of $\rank r$ and $\deg d$
where $(r,d)=1$ form a fine moduli space.
In particular, there is a Poincare bundle $\Cal P$
over the moduli spaces $SU$ and $U$.  
In the non-coprime case, the variety is quasi-projective 
and has no Poincare bundle (even locally).  
However $U$ and $SU$ can be covered by Zariski
open sets $\{A_i\}$.  The $A_i$ in turn have etale
covers $p_i: T_i \to A_i$. There is a Poincare bundle
$\Cal P_i$ over $T_i$ such that the map given by
the universal property of $U$ and $SU$ is $p_i$.
This gives a parameter space of stable bundles
with rank $r$ and degree $d$ where $(r,d)\ne 1$.

Sundaram constructs parameter spaces from the Jordan
Holder filtration, and the Harder Narasimhan filtration.
The strictly semistable bundles are  parameterized
by finitely many components with dimension $\le \dim M -(r-1)$
The unstable bundles have infinitely many components.
So this case needs explaining.  First of all, for any
Bundle E with rank r and degree d, there is a unique
filtration:
$$ \sum E: 0 = \Cal E_1 \subset \Cal E_2 \dots \subset \Cal E_k = E,$$
where $E_i = \Cal E_i/\Cal E_{i-1}$ is semistable, and 
$\mu_i =\mu(E_i)$ is strictly decreasing.
We say that $E$ has Harder Narasimhan (or HN) type
$\{(r_1, \mu_1), (r_2, \mu_2), \dots (r_k, \mu_k)\}$;
$r_1+r_2+ \dots +r_k =r$, where $r_i=\operatorname{rank}(E_i)$.  
We set $\mu_1 = \mup$;
it is the upper bound on $\mu$ of a subbundle.
Similarly, we set $\mu_k = \mdn$; it is the lower
bound on $\mu$ of a quotient bundle.
If $\mup < 0$ then $h^0=0$ as $\Cal O_C$ is
not a subbundle.  And by Serre Duality,
if $\mdn(E)>2g-2$, then $h^1=0$.
Since the complex numbers are characteristic zero,
$\mup(E\otimes F)=\mup(E)+\mup(F)$.
$\mdn(E\otimes F)=\mdn(E)+\mdn(F)$.
Also $\mup(E^*)=-\mdn(E)$ by dualizing the HN filtration.
and likewise $\mdn(E^*)=-\mup(E)$.

If $r$ is fixed and $\mup$ is bounded from above
and $\mdn$ is bounded from below, there are only
finitely many HN types possible.  
This is also the case when $\deg$ is bounded from
above and $\mdn$ is bounded from below, or if
$\mup$ is bounded from above and $\deg$ is bounded
from below.  The advantage is that then there are 
only finitely many  components of the parameter space
of unstable bundles, and each has dimension $\le \dim U -(r-1)$
or $\le \dim SU -(r-1)$.

\heading  \S 2 Some Propositions \endheading

\proclaim {Proposition 1}
Let $0 < c \le d$ and $c<r$ and $L$ be a line bundle of $\deg =0$.  
For a general $E_t\in U(r,r(g-1)+d)$ with $d>0$ and $d$ not
necessarily $< r$, there is an exact sequence:
$$0 \to \overset c \to \oplus L \to E_t \to F_u \to 0$$
with $F_u$ stable.
\endproclaim

\subheading {Remark 2}  This does not mean that all embeddings
of $\cplus L$ have a stable  vector bundle cokernel.  
But a general subbundle of a general bundle such a cokernel.  
Of course if $c=d$ there is only one such subbundle 
for general $E_t$ and the cokernel is unique.

\demo{Proof of Proposition 1}
The first thing that can go wrong is that the cokernel
is not a bundle.  In some sense this is all
that can go wrong outside of stability.
If $c=1$ this is obvious.  If $c=2$,
the one dimensional family of (twisted) sections may
span a subbundle of $E_t$ of rank 1.  But going back to
$c=1$, that means $L$ factored.  If that never happened,
there is no problem at $c=2$ except maybe factoring through
a subbundle of rank $2$, and so on. 
So we can assume $\cplus L$ is a subbundle.

So now the idea is to count dimensions of the space of
subbundles in various $E_t$ that our trivial bundles
could factor through and the dimensions of the space of 
subbundles of the form $\cplus L$.
We know $\cplus L$ embeds in every $E_t$.
We also know $h^0(C, E_t) =d$ for general $t$
and there are $c(d-c)$ dimensions of embeddings.
So the total dimension of subbundles of the form $\cplus L$ is 
$$r^2(g-1) + 1 + c(d-c), \tag $*$ $$ 
which we will show is larger than the dimension 
of the space of bundles factored through.

So suppose $\cplus L$ factors through a bundle
$S$ and that we have a sequence;
$$0\to S \to E_t  \to  Q \to 0 $$
where $Q$ is a vector bundle.  $S$ is a rank $c$
bundle and an elementary  transformation of $\cplus  L$.
So at  first glance the dimension of possible bundles $S$
is $c\delta$ with $\delta = \deg S$, since $S$ must be
an elementary transform of $\cplus L$.   But there are $c^2-1$
projective automorphisms of $\cplus L$.  
So if  $\delta \ge c$ the dimension is  $c\delta - c^2+1$.
If $\delta < c$ the dimension is $0$ and there are at least
$c^2-c\delta-1$  dimension of automorphism of $S$.  The  dimension
of automorphisms  must be subtracted from the dimension
of extensions, so we shall abuse terminology and say the
dimension of the space of bundles $S$  with a map 
$\cplus L \to S$  is $c\delta -c^2+1$.

We know that 
$$\deg Q = r(g-1) + d -\delta = (r-c)(g-1) + c(g-1) +d -\delta.$$
And since $\mudn (Q) \ge \mu (E_t)>0 $ 
is bounded from below, only finitely
many HN types are possible for $Q$.  The dimension is bounded
by the dimension for stable  $Q$, that is $(r-c)^2(g-1) +1$.
Since $Q$ is a quotient of the stable bundle $E_t$.
$\mdn(Q)>\mu(E_t)$.  Similarly $\mup (S)<\mu(E_t)$
And hence $\mudn(S^*)>-\mu(E_t)$.  All of which gives
us $\mudn (\omega_C \otimes S^*\otimes Q) > 2g-2$
And hence, $\hi \omega_C\otimes S^* \otimes Q)=0$.

Now we estimate and bound the dimension of the space of sections.
$$ h^0(C, \omega_C \otimes S^* \otimes Q)  =
2c(r-c)(g-1) -\delta(r-c)+ c^2(g-1) + cd -c\delta.$$

Now  we add the dimension of bundles $S$ and $Q$
then subtract $1$ from $\ext$ in order to
get the projective dimension.   It all comes out to:
$$r^2(g-1) +1 + cd -c^2 -\delta(r-c) .\tag ${*}{*}$ $$
The upshot is that $*$ is bigger than ${*}{*}$,
so $\cplus L$ is a subbundle with a bundle cokernel.

To see that the general cokernel is stable, consider
yet another dimension count.  We  count extensions of the form
$$ 0 \to \cplus L \to E_t \to F_u \to 0.$$
Let $F_u\in U$ be a parameter space of stable bundles of
$\rank r-c$ and $\deg r(g-1) + d$.  
It has dimension $(r-c)^2(g-1) +1$.
Calculating as above we see the dimension of extensions
is $r^2(g-1) +c(d-c)$.  A  parameter space of nonstable 
bundles would have dimension at least $n-1$ lower;
therefore, the general cokernel is stable.
That proves Proposition 1. \qed
\enddemo

\proclaim {Proposition 2}  Suppose that for a general bundle
$E_t\in U$ there is a natural sequence
$$0\to W_k \to E_t \to F_u \to 0$$
where $F_u$ is a stable vector bundle
for general $u$.  If $\cplus L$ (with $\deg L =0$)
embeds in a general $E_t$,
$L \nsubseteq W_k$, 
and $dk + c < r$,
then for (perhaps more general $E_t$) there is an exact sequence:
$$0\to W_k \cplus L \to E_t \to F^{'}_v \to 0$$
with $F^{'}_v$ a stable vector bundle
and $\cplus L$ a general subbundle of $E_t$.
Hence, by setting $c=d$ and considering only 
$E_t$ with $h^0(C, E_t \otimes L_i^*)=d$,
and with $ed < r$, we get for general t, 
a unique sequence: 
$$0\to W_e \to E_t \to F_u \to 0$$
with $F_u$ a stable bundle for general $t$.
\endproclaim

\demo{Proof of Proposition 2}
First we start with $F_u=F$ fixed.
Let $\cplus L$ be a subbundle of $F$.
Then for a general subbundle Proposition 1 applies
and the cokernel is stable.  Now we have to
pullback the subbundle of $F$ to a subbundle
of $E_t$.  

For any given embedding $\cplus L @> \alpha >> F$,
there is a commutative diagram:
$$\CD
0 @>>> W_k @>>> E_{t\alpha} @>>> \cplus L        @>>> 0 \\
@.    @VVV      @VVV        @VV\alpha V          @.  \\
0 @>>> W_k @>>> E_t   @>>>     F                 @>>> 0 
\endCD 
$$
$E_{t\alpha}$ is base extension of $E_t$ over $F$ by $\alpha$.
So $\alpha$ factors through $E_t$ iff it factors through $E_{t\alpha}$,
which means the top sequence splits.
To interpret this, consider the surjective induced map:
$$H^0(C, \omega \otimes W_k^* \otimes F)^* \twoheadrightarrow 
H^0(C, \omega \otimes W_k^* \otimes (\cplus L))^*.$$
The top sequence splits if the element corresponding
to the bottom sequence maps to $0$.
It is easily verified by a dimension count that there
are extensions which cause the upper sequence to split,
and many of them.  For the kernel is the dimension of
the bottom extensions of $F$ to $W_k$
minus the dimension of top extensions $\cplus L$ to $W_k$.

We can choose $E_t$ so it has  
$\hi E_t \otimes L_i^*)=0$ for $1 \le i \le k$ 
and $\hi E_t \otimes L^*)=0$.
This forces a choice of $F$ but we can choose $E_t$
so that $F$ is stable.  Also note, the general extension
of $F$ by $W_k$ will be stable and have the above vanishings
of $h^1$.
Now we count the dimension of pairs,
$(\cplus L \subset E_t)$ denoting a specific subbundle
$\cplus L$. It's dimension is the dimension of the bottom extensions
plus $c(d-c)$.  Now this maps to the set of subbundles
of $F$.  So we subtract the dimension of the kernel,
which means we subtract the dimension of the bottom extensions
and add the number of top extensions.
This gives us a total of $c(d-c) + kcd(g-1)$.

Now we count the dimension of pairs $(\cplus L \subset F)$.
Since a general $E_t$ has $h^1(C, E_t \otimes L^*)=0$,
so does $F$ have $h^1(C, F \otimes L^*)=0$.
Since $h^0(C, F \otimes L^*) = r(g-1) + d - (r-dk)(g-1)$,
we get dimension of rank $c$ subbundles given by,
$cdk(g-1) + c(d-c)$.  So the map is dominant and a general
bundle $\cplus L \subset E_t$ has a stable cokernel
(provided our bundle is an extension of $F$).

The preceding argument has two faults.
The first is that we never constructed the
object whose dimensions we counted.
The second is that it handles $F_u$ one $u$ at a time,
treating the family of $F_u$ as a set 
(with no algebraic structure).

For the first, let $V= \ext (F, W_k)$.
Now there is a sequence over $V\times C$:
$$0\to \Cal W_k \to \Cal E \to \Cal F \to 0.$$
Given a point $v\in V$ the sequence restricts
over $v\times C$ to:
$$0\to W_k \to E_t \to F \to 0,$$
where the extension corresponds to $v\in V$.
See \cite{5} \cite{7} and \cite{4}.
If we tensor by $L^*$ and push down 
$\pi: V\times C \to V$, we get a sequence
$$0 \to \pi_* \Cal E^{'} \to \pi_* \Cal F^{'}.$$
When restricted to $v\times C$ we get:
$$0\to \bho E_t) \to \bho F).$$
A rank $c$ subspace in $\bho E_t)$ corresponds
to a unique subbundle  $\cplus L \subset E_t$.
And its image in $\bho F)$ corresponds to
the subbundle's image in $F$.
So we take the Grassmannian bundles and note
that the second Grassmannian bundle maps from
$V\times Gr$ to $Gr$.  So a Grassmannian bundle
over $V$ dominates a Grassmannian.

The previous dimension count now holds.
And our dimension count has shown a birational
map between subbundles of general $E_t$
and those of $F$.  We have a general
subbundle, $\cplus L \subset F$ has a stable cokernel by Proposition 1.
And that implies a general subbundle $W_k \cplus L \subset E_t$
has a stable cokernel provided $E_t$ has cokernel
$F$ corresponding to $W_k \subset E_t$.

We need to extend this result to a general cokernel $F_u$.
So we need a general construction of a universal
extension of two families of vector bundles.
We have $W_k$ which does not really vary,
and $F_u$ which does. So we have two families $T$
and $S$.  $T$ is a point, and $S$ is some cover
of an open subset of $SU$ the moduli of stable
bundles of some fixed rank and degree.
First we need $\bhi F_u \otimes L^*)=0$.
This is true for general $F_u$.
Then we need to take an open set and an etale
cover to get a Poincare bundle.  We then take
$S$ a Zariski open affine subspace of the etale
cover.  This will still dominate.
We now satisfy the criteria of \cite{Lemma 2.4, 7},
because an affine space has no higher cohomology.
So by that lemma, there is a universal extension
with parameter $\prj(V)$ which is a projective bundle over $U$
whose fiber at $u$ is $\bhi \operatorname{Hom}(F_u, W_k) ) $.
It can be formed by taking the pullback of $\omega_C\otimes W^*$,
tensoring by the Poincare bundle $\Cal P$ and pushing down
onto $U$ and dualizing.
Now we can take Grassmannian bundles and globalize the
proof, showing that for a general subbundle
$\cplus L \subset E_t$ of a general bundle $E_t$
the cokernel is a stable bundle. 
\qed
\enddemo

\proclaim{Proposition 3}
For a general bundle $E_t\in U$, 
a line bundle $L_{k+1}$ of $\deg 0$ 
with $L_{k+1} \nsubseteq W$, 
a general $\cplus L_{k+1} \subset E_t$ with $c\le d$,
and $W_k \subset E_t$, the natural map
$$W_k \cplus L_{k+1} \to E_t$$ is surjective
provided $dk + c \ge r+1$.
\endproclaim

\demo{Proof of Proposition 3}
To begin with we show that the map is at least generically
surjective for a general $E_t$.
So let $i$ be the largest integer such that $di \le r-1$.
Then for a general $E_t$ and a general $\cplus L_{i+1} \subset E_t$,
with $c^{'}=r-1-di$, 
we get a subbundle whose quotient is a line bundle,
by Proposition 2.
Now another subbundle $L_{i+1} \subset E_t$ which does
not factor through $W_i \overset c^{'} \to \oplus L$
maps to the line bundle quotient generically.
We can choose another $L_{i+2}$, set $c^{'} + 2 = c$,
and the map $W_k\cplus L_{k+1}$ generically surjects
onto $E_t$.

Now we have to show surjectivity for a general $E_t$.
Once again, this will be done by dimension count.
First we count the number of pairs $W_k \cplus L_{k+1} \subset E_t$.
The dimension of the $E_t$ is $r^2(g-1) +1 =(r^2-1)(g-1) + g$.
The dimension of the subbundle (for given $E$) is the dimension
of the $\operatorname{Grassmannian}(c,d)=cd -c^2$.  So our final
dimension is:
$$(r^2 -1)(g-1) +g +cd -c^2. \tag{${*}{*}{*}$} $$

Now suppose the map in the proposition does not surject
onto $E_t$.  It then surjects onto some bundle $Q$ and
$E_t$ is an elementary transformation of $Q$.  So
this means 
$$\deg Q = \deg E_t - \delta = r(g-1) +d -\delta $$
To calculate the dimension of the family of $Q$'s
consider the sequence:
$$0 \to K \to W_k \cplus L_{k+1} \to Q \to 0,$$
where $K$ is a line bundle, with $\deg K = r(g-1) +d -\delta$.
Assume for now $\hi K^*L_i)=0$ for all $L_i$ and $K$.
An informal dimension count gives:
$$(r+1)(r-1)(g-1) +g +(r+1)d - (r+1)\delta.$$
This follows from there being $r+1$ subline bundles in
$W_k \cplus L_{k+1}$.  All of them have $\deg = 0$,
and $K^*$ has $\deg = r(g-1) +d -\delta$, and therefore,
$\ho K^*\otimes L_i)= (r-1)(g-1) +d -\delta$.
The $g$ is the dimension of the possible $K$.
Now we need to account for the automorphisms
which come to $kd^2 -c^2$.  Now note that 
$(r+1)d = (kd +c)d = (kd^2 +cd)$.  And hence,
$(r+1)d-kd^2 -c^2 = cd -c^2$.  Finally we have
to add $r \delta$  the dimension of elementary
transformations.  The total dimension is:
$$(r^2-1)(g-1) +g +cd -c^2 -\delta.  \tag ${*}{*}{*}{*}$.$$

Since $({*}{*}{*}) \ge ({*}{*}{*}{*})$ with equality iff $\delta = 0$,
we can conclude surjectivity except for two problems.
The first: what object have we calculated the dimension of?
The second: what if $K^*\otimes L_i$ is special?

For the first problem we can form a universal
Hom space, for maps from $K$ to $W_k \cplus L_{k+1}$.
The quotient comes for free.
In general the Hom space is constructed as follows.
Suppose we have families of bundles $\Cal E$ over
$S\times C$, and $\Cal F$ over $T\times C$. 
Now consider $\Cal E^*\otimes \Cal F$ over 
$S\times T \times C$, and push down from $C$,
then pull it back up and call it $U$.  
If $h^0(C, E^*\otimes F)$
is constant, then restricting to a point $(E\in S, F \in T)$ 
the bundle is $\bho E^*\otimes F)$.  And over a point
$v\in \bho E^*\otimes F)$ we get a morphism $E \to F$
over $C$.  
If $h^0(C,E^*\otimes F)$ is not constant we can take
a stratification and check the dimension of each piece.

For the second problem we have $\deg (K^*)=\alpha$
with $1 \le  \alpha \le 2g-2$.  So $g\ge 2$.
We will do the cases $\alpha = 1$ and $2$ 
separately.  Then we will do the case $\alpha \ge 3$.
For $\alpha =1$, $d$ must be $1$ and $k\ge 2$.
For fixed $L_1$ and $L_2$ distinct there are only
finitely many line bundles $K$ of $\deg = 0$
so that $\ho K^*\otimes L_i)\ne 0$ for $i=1$ and $2$.
Furthermore we must have $d=1$.
The morphisms from $K$ to $L_i$ are unique up to
automorphisms of $W$ and hence there are only finitely many $Q$,
one for each $K$.
Now the dimension of elementary translations
is $r$ times the difference of degree or
$$ r(r(g-1) +d -1) <  {*}{*}{*}.$$

Now we assume $\alpha =2$ and $h^0(C, K^*\otimes L_i)=1$ for all $i$,
and hence $d=1$.
The dimension of $K$ is $\le \alpha$.
As above each $K$ gives rise to a unique $Q$.
So the dimension of $Q$ is $\le \alpha$ and
adding the dimension of elementary transformations we get
$$\alpha + r(r(g-1) +d-\alpha) < {*}{*}{*}.$$

Now we by reording $i$ if necessary, we can assume
$h^0(C,K^*\otimes L_1)\ge 2$  and $\alpha =2$.
$h^0(C,K^*\otimes L_i)=1$ for $i\ne 1$.
We have two case.  Case i) d=2 and hence k=1.
Case ii) d=1.  For case i, the dimension of $K$ is zero.
Furthermore, the dimension of $Q$ is zero because there
is only one fixed map from $K$ to $W$ up to automorphisms of $W$. 
Finally, the dimension of elementary transformations is 
$$r(r(g-1)+d-2) < {*}{*}{*}.$$
For case ii, dimension $K$ is again zero,  but there
is a one dimensional family of maps from $K$ to $L_1$ 
up to automorphsims of $W$,
and hence, there is a one dimensional family of $Q$.
We add the dimension of elementary transformations
to get
$$1+r(r(g-1)+d-2) < {*}{*}{*}.$$

\comment
Now assume $\alpha = 3$ and hence $g\ge 3$ as
$\alpha \le 2g-2$.  There are two cases:
i) $d= 2$, and ii) $d=1$.  In case i) we have
a 1 dimension family of bundles $K$ now if $c=2$
there is for each $K$ one $Q$.  But if $c=1$ there
is for each $K$ a 1 dimension family of $Q$.
So dimension of $Q$ is dimension of $K+ cd -c^2$.
Adding the dimension of $Q$ to the dimension of
elementary transformations we get:
$$1+r(r(g-1)+d-3) +cd-c^2 < {*}{*}{*}.$$
For case ii, there is a 1 dimensional family of $K$
and $r+1$ distinct bundles.  There is a 1 dimension
family of maps from $K$ to $L_i$ for each $i$.
So the dimension of $Q$ is $r+2$.  Keeping in mind
$d=1$ we add the dimension of $Q$ to the dimension of
elementary transformations and we get:
$$(r+2) +r(r(g-1) + d -3) < {*}{*}{*}.$$
\endcomment

Now assume $3\le \alpha \le 2g-2$.
By Clifford's Theorem, 
$\ho K^*\otimes L_i)\le \frac{\alpha}{2} +1$.
By Marten's Theorem the dimension of $K$ with
$\ho K^*\otimes L_i)\le \frac{\alpha}{2} +1-a$
equals $2a$ for $a$ non-negative.  Now assume all $L_i$ have
$\ho K^*\otimes L_i)\le \frac{\alpha}{2} +1-a$.
The dimension of $K$ is $2a$.  The dimension of $Q$
is $\le kd(\alph +1 -a -d) + c(\alph +1 -a-c) + 2a$ or:
$$(kd+c)(\alph +1 -a -d) +cd - c^2 + 2a.$$
This is clearly maximal if $a=0$, as $kd+c=r+1>2$.
Now adding the dimension of the elementary transformations
and rearranging the terms we get:
$$r^2(g-1) + cd - c^2 +\alph(r+1)-r\alpha -d +r+1.$$
To show this is less than ${*}{*}{*}$ it suffices to show:
$$
\align
\alph(r+1)-r\alpha -d +r+1        &< 1 \\
\frac{-(r-1)}{2}\alpha -d +r+1    &< 1 \\
-\frac{3}{2}r+\frac{3}{2}-d+r+1   &< 1 \\
\frac{-r}{2} -d +\frac{5}{2}      &\le \frac{1}{2} \qed
\endalign
$$
\enddemo

\heading \S 3 The Main Theorem \endheading

Now we are prepared to prove the Main Theorem.

\demo{Proof of Main Theorem part A}
Suppose $ed = r-1$.  By Proposition 1 and repeated use of
Proposition 2, for a general $E_t \in SU$ we have a sequence:
$$0\to W \to E_t \to A \to 0$$ for a fixed line bundle $A$.
Now the family of $E_t$ is given by elements of 
$\ext (W, A)$ modulo the automorphism of $W$.       
$\ext$ is given by $\bho \omega_C \otimes W^* \otimes L)^*$.
This decomposes into a direct sum of 
$\bho \omega_C \otimes V_i^* \otimes A)^*$,
which in turn is isomorphic to a direct sum of
$\bho \omega_C \otimes A \otimes L_i)^*\otimes U_i^*$.
Concentrating on one $i$ we can see that the
quotient by the automorphisms of $U_i$ gives a Grassmannian.
So the whole thing is just a product of Grassmannians.

Now for $ed = r+1$ , Proposition 3 shows that for a 
general $E_t \in SU$, we have an exact sequence
$$0\to A \to W \to E_t \to 0,$$
for some fixed bundle $A$.
So the family is parameterized by
$\operatorname{Hom} (A, W)$ modulo the automorphisms
on $W$. We can decompose and get a direct sum
of $\bho L^*\otimes L_i)\otimes U_i$ modulo the automorphisms of $W$ ---
but again, this is just a product of Grassmannians.  
And we are done.
\qed 
\enddemo

To prove part B we need a result of Boden and Yokogawa.

\subheading{Definition}  A variety $X$ is stably rational
of level $k$ if $X\times \Bbb P^k$ is rational.

\proclaim {Theorem 3 (Boden and Yokogawa \cite {3}) }
If $(r,d)=1$ then $SU(r,d)$ is stably rational
of level $r-1$
\endproclaim 

Now we prove part B.

\demo{Proof of Main Theorem part B}
As in the argument for part A, using Propositions
1 and 2 repeatedly we get a sequence for general
$E_t$,
$$0 \to W \to E_t \to F_u \to 0,$$
where $F_u$ is stable.  This means a general
$E_t$ is an extension of a stable bundle $F_u$
by $W$.  So for a general $F_u$ and a general
extension by $W$ we get $E_t\in SU$.  The map
to $SU$ is not one to one, but its fiber
is just the automorphisms of $W$ since for general
$E_t$, $W$ is a unique subbundle.

So now we need to construct the map more carefully,
then take care of the fiber.  First note that 
by the hypothesis the rank and degree of $F_u$
(which is $(\delta, r(g-1) +d)$) are coprime
because all factors of $\delta$ divide $r$ and not $d$.
So now let $SU^{'}=SU(\delta, r(g-1) +d)$.
$SU^{'}$ is a fine moduli space with a Poincare
bundle $\Cal P$ over $SU^{'}\times C$.
We also have projections $p_1$ and $p_2$ onto 
$SU^{'}$ and $C$ respectively.  Now we will form
a bundle parameterizing extensions of $F_u \in SU^{'}$ by $W$.  
Consider 
$$p_2^*p_2{}_*(p_1^*(W^*\otimes \omega_C) \otimes \Cal P)^*.$$
The general fiber is $\bho W^*\otimes \omega_C \otimes F_u)^*$
which is canonically $\ext (F_u, W)$.
It is also canonically 
$\overset e \to \oplus \bho L_i^* \otimes \omega_C \otimes F_u)^*
\otimes U_i^*$.  By taking an affine open subset on $SU^{'}$,
we can assume the bundle is trivial, and that we have a universal
extension space parameterized by the total space of the bundle,
which by restricting to a smaller affine open subspace of $SU^{'}$
is trivial and hence has no higher cohomology which satisfies
the conditions of \cite {Lemma 2.4, 7}.  And therefore
this gives a dominant map from the bundle over $SU^{'}$ to $SU$
which has fiber isomorphic to $\aut (W)$.  But clearly the
action of $\aut (W)$ on the bundle over $SU$ gives a product of
Grassmannians cross $SU^{'}$.  The dimension of the product
of Grassmannians is the rank of the bundle minus $ed^2$.
For the rank of the bundle we have 
$$h^0(\omega_C \otimes L_i^* \otimes F_u) =
\delta(g-1) + (r-\delta)(g-1) + d +2(g-1)\delta -\delta(g-1).$$
This is $(r+\delta)(g-1) + d$.  Now to get the rank of the full
bundle we multiply by $e$ the number of $L_i$, and $d$ the rank
of $U_i$.  But $ed=r-\delta$, so we get
$(r+\delta)(r-\delta)(g-1) +ed^2$.  Now $\aut(W)$ has dimension $ed^2$.
So the product of Grassmannians has dimension
$(r+\delta)(r-\delta)(g-1)\ge r-1$.
Now by Boden Yokogawa $SU^{'}$ cross the product of Grassmannians
is rational.  And therefore so is $SU$.
\qed
\enddemo

\subheading{Remark 3}  The condition that $ed=r-\delta$
where all prime divisors of $\delta$ are divisors of
$r$ is necessary.  Suppose to the contrary that a prime
$p$ divides $\delta$ but $p$ does not divide $r$, 
then we can solve the equation $r(g-1)+d \equiv 0 \pmod p$.  
The rank and degree would then not be relatively prime for
those $g$ which solve the equation.  

\subheading{Remark 4}  Solutions for $g$ that satisfy the
above remark are the only $g$ where rationality can fail,
because otherwise $(r(g-1)+d, \delta)=1$ and the proof
goes through.  For a given $g$, it is enough if 
$(r(g-1)+d, \delta)=1$ for one possible $\delta$.
In the case of $r=13$ and $d=5$ or $d=8=(13-5)$
(the first case where $(r,d)=1$ but the hypothesis
of the Main Theorem does not hold) we get three 
equations.  First the Boden Yokogawa equations
$(g,5)=1$ and $(g,8)=1$.  The Chinese Remainder
Theorem gives $g\equiv 0 \pmod {10}$.  Now we
let $e=2$ so $\delta = 13 - 2(5)=3$ and we
get $13(g-1)+5\equiv 0 \pmod 3$.  The final equation is
then
$$
g\equiv 20 \pmod {30}.
$$
We shall have more to say about this in the next section,
after making a conjecture.

\heading \S 4 Conjectures \endheading

Inspired by the Main Theorem, it seems natural
to prove that $SU$ is rational if $ed = r +\delta$
where $\delta$ has the property that any prime
divisor of $\delta$ is a prime divisor of $r$,
and $0<d<r$.

We {\bf attempt \/} a proof.

We begin by constructing the universal hom space
over $SU$  and \{W\}.  It might have more than 
one irreducible component, corresponding to
bundles $E_t$ where $\ho W^*\otimes E_t)$
is large because $h^1$ is large.  But for a general bundle we have
the above $h^1=0$.  So we get one irreducible
component which dominates $SU$.
Now given a morphism $W\twoheadrightarrow E$.
We get a morphism $K_{E,W}\hookrightarrow W$.
We would like $K_{E,W}$ to be stable.  And given
$SU^{'}=SU(\delta, -r(g-1)-d)$, there is a dominant
map from $SU$ to $SU^{'}$.  This would imply that
$\delta < r$ since if it is bigger than $r$
the map $SU$ to $SU^{'}$ cannot dominate for dimension
reasons.  And if $\delta = r$ we get birationally isomorphic
spaces $SU$ and $SU^{'}$. 

But how do we know $K_{E,W}$ is stable for a general
$E$ when $0< \delta <r$?  The problem is we do not!  
So for now, {\bf suppose \/} that $K=K_{E,W}$ is stable.
For a general $K$ we have $\hi K^*\otimes W)=0$.
So we calculate the dimension of the irreducible component
of the universal Hom space containing $\{K | \hi K^*\otimes W)=0\}$.
First we count the dimension of $\operatorname{Hom} (K,W)$.
$$
\align
\ho K^*\otimes W) &= \delta ed(\frac{r(g-1)+d}{\delta} - (g-1)) \\
                  &= ed(r-\delta)(g-1) + ed^2 \\
                  &= (r^2-\delta^2)(g-1) + ed^2
\endalign
$$
Adding $\dim U^{'} = \dim U(\delta, -r(g-1)-d) =(\delta^2(g-1))$, we get
$r^2(g-1)+ed^2$.  This is the dimension of $U$ plus
the dimension of $\aut W$.  So we get a dominant
map from the irreducible component of the Universal Hom
space dominating $SU{'}$ to $SU$ and the general fiber is
$\aut(W)$.  

Since $(\delta, -r(g-1)-d) =1$ because
all prime divisors of $\delta$ divide $r$,
we can then restrict to an open set of $SU^{'}$ 
and construct the Universal Hom bundle by
Taking the pullback of the pushdown of the Poincare bundle $\Cal P$
tensored by the pull back of $W^*$.  This bundle
is locally trivial and has fibers $\bho E \otimes W^*)$
which equals the direct sum $\overset e \to \oplus
\bho E \otimes L_i^*)\otimes U_i$.  
Now the action of $\aut (W)$ gives $SU^{'}$ cross a product
of Grassmannians. Now $SU$ is rational because of Boden Yokogawa
and the fact that the dimension of the product of Grassmannians is 
$\ge r-1$.

\proclaim {Conjecture 1}  Let $C$ be a smooth projective irreducible
curve over $\Bbb C$.  If $(r,d)=1$ and $d$ divides $r+\delta$
where $0<\delta <r$ and every prime divisor of $\delta$ divides $r$,
then $SU=SU(r, r(g-1)+d)$ is rational.
\endproclaim

As noted in the above plausibility argument, 
Conjecture 1 would follow from the 
stability of the bundle $K_{E,W}$ and the vanishing of 
$\hi K_{E,W}^*\otimes W)$ 

But it is actually enough to show the vanishing 
of one $\hi K_{E,W}^*\otimes W)$ provided $\hi W^*\otimes E)=0$.  
This is because we would then get a family of $K_{E,W}$
for which  $\hi K_{E,W}^*\otimes W)=0$ generically.  The dimension
of that component of the Universal Hom is, by the calculations
preceding Conjecture 1, the dimension of the family
$+(r^2-\delta^2)(g-1) + ed^2$.  Now this adds up to 
$r^2(g-1) +1 + ed^2$ (the dimension of $\operatorname{Hom}(W,E)$)
if and only if the dimension of the family is $\delta^2(g-1) +1$.
But $\deg K_{W,E}$ is fixed, and $\mup<0$ because it is a
subbundle of $W$.  This means there are only finitely many
HN types.  But the dimension of a parameter space for unstable bundles
is less than $\delta^2(g-1) +1$.
Hence, $K_{W,E}$ is generically stable.
So it would suffice to prove the following conjecture
in order to prove Conjecture 1.

\proclaim {Conjecture 2}
Let $C$ be as above.  There is an $E\in SU$
such that $\hi W^*\otimes E)=0$, $E$ is spanned by $W$,
and $\hi K_{W,E}^*\otimes W)=0$.
\endproclaim

We can prove Conjecture 2 if $e=1$ or $2$.

\proclaim {Proposition 4}  Conjecture 2 holds if
$e=1$ or $2$.
\endproclaim

To prove this we need a Lemma belonging to folklore.
By an elementary transformation of a bundle $E$ we mean
the kernel, $E_p$, of a sequence. 
$$0\to E_p \to E \to \Cal O_p \to 0,$$
where $p$ is a reduced point on $C$.
This corresponds to blowing up a point $p\in \Bbb P(E)$
and taking the pushdown of the fiber.
Corresponding to a point on $\Bbb P(E)$ there
is a $r-2$ plane in $\Bbb P(E_p)$ and this corresponds
to a point in $\Bbb P(E_p^*)$.  If we blow up this
point, blow down the fiber, and then dualize we get
$E$ back.  This is known as an inverse elementary
transformation.

\proclaim {Lemma 1} Let $C$ be  as above.   For a general
bundle $E\in U(r,d)$ or $E\in SU(r,d)$ a general elementary
transformation is stable and so is a general inverse elementary
transformation.
\endproclaim

\demo {Proof of Lemma 1}  Suppose to  the contrary that
the general $E_p$ and hence all $E_p$ are unstable.
$E_p$ has fixed degree $(d-1)$ and $\mup (E_p) < \frac{d}{r}$
because $E_p$ is a  subbundle of the stable bundle $E$.
This means there are finitely many HN types.

Now consider the parameter space of a general $E_p$
with it's universal bundle $\Cal E$.  Take $X=\Bbb P(\Cal E)$.
A point on $X$ is an elementary transformation of a bundle $F^*$
which when dualized gives a bundle $E\in U(r,d)$.  If the map 
to $U(r,d)$ does not dominate then the Lemma is proved.  So assume
the map dominates.  The dimension of $X$ is the dimension of
the parameter space $+ r-1$.   Which equals (by Sundaram)
$r^2(g-1) -r+1$ for the parameter space $+r-1 = r^2(g-1)$.
That means the fiber for a general bundle $E\in U(r,d)$
consists of finitely many points.  But those points represent
all the elementary transformations of $E$ which are not stable
(because every such transformation has an inverse).
The conclusion is a general elementary transformation of
a stable bundle is stable.
This proves the Lemma.
\qed
\enddemo

\demo {Proof of Proposition 4}
In the case $e=1$ there is only one linebundle $L_1$.
Twisting by $L_1^*$ we can assume the bundle is $\Cal O_C$
and $W$ is $H^0(C, E)\otimes \oc$.
We will now construct a bundle for $d= 2r$
then work our way down in degree.
Let $A$ be a general line bundle of $\deg g+1$.
$A$ is non special and spanned.  So let $E$
be the direct sum of $r$ copies of $A$.
$K_{E,W}^*=E$ and so it is nonspecial.
The only problem is that the original bundle $E$
is not stable.  We can choose a one parameter
family of stable bundles $E_t$ degenerating to $E_0=E$.
A general $E_t$ will be stable and spanned and so
will have a bundle $K_{{E_t},W}$.  So taking
a possibly smaller family we parameterize a family
of $K_{{E_t},W}^*$ with a general bundle nonspecial.
This finishes the case $d=2r$.

To get the case $d=2r-1$ we use elementary transformations.
We can assume that a general $E_t \in M_r=U(r,rg+r)$ has
$K_{E,W}^*$ non special.  Furthermore, a general elementary
transform of $E_t$ is stable by Lemma 1.  And for $F\in M_{r-1}=U(r,rg+r-1)$,
a general inverse transformation is stable by Lemma 1, and a general
bundle $F$ is spanned.  Putting these together shows that there
is a pair $E \in M_r$ and $E_p \in M_{r-1}$.  Where $E$ is spanned
with $K_{E,W}^*$ nonspecial and $E_p$ an elementary transformation
of $E$ with $E_p$ spanned.
We have an exact sequence
$$0\to E_p \to E \to \Cal O_p \to 0.$$
There is a second exact sequence
$$0 \to \bho E_p)\otimes \oc \to \bho E)\otimes \oc \to \oc \to 0.$$
The second sequence surjects onto the first and the kernel is
$$0 \to K_{E_p,W_p} \to  K_{E,W} \to \oc(-p) \to 0.$$
$W_p$ is just $\bho E_p)\otimes \oc$.  And by dualizing the
last sequence and taking cohomology we see that 
$K_{E_p,W_p}^*$ is nonspecial.  Repeating this argument
proves the case $e=1$.

For the case $e=2$ there are line bundles $L_1$ and $L_2$ of
$\deg 0$.   Now do the case where $d=r$.  Choose a general
line bundle $A$ with $\deg(A)=g$ so that $A\otimes L_i^*$
is not special for $i= 1$ or $2$.  And choose $A$ such that
the divisors $|A\otimes L_1^*|$ and $|A\otimes L_2^*|$
have no points in common.  Now let $E$ be the direct sum
of $r$ copies of $A$.  It is easily verified that 
$\hi K_{E,W}^*\otimes L_i)=0$ for $i = 1$ or $2$.
This proves the case $d=r$.  The remaining cases are
proved as above using elementary transformations.
\qed
\enddemo 

\subheading{Remark 5}
The case $e=1$ has no impact on Conjecture 1.
The case $e=2$ seems to but does not.
For suppose that $2d = r +\delta$ where every
prime factor of $\delta$ divides $r$.  Then 
$2(r-d) = r - \delta$  where $\delta$ has the
above properties.  So the result is already
known because $SU(r, r(g-1)+d)\cong SU(r, r(g-1) +r-d)$.

However the case $e=1$ may have some use.
Throughout this paper we assumed the choice of bundles
$\{L_1, \dots L_e\}$ was of arbitrary distinct bundles.
But we never needed that fact (outside of distinctness).
So consider $L_1$ such that $L_1^{\otimes i}$ is trivial
iff $i$ is a multiple of $e$.  Now set $L_i=L_1^{\otimes i}$.
The set of $L_i$ corresponds to a cyclic \'etale cover
$\tilde C$ of $C$.  The $L_i$  pull back to be trivial,
so $\tilde W$ the pullback of $W$ is trivial.
On $\tilde C$, a general $K_{\tilde E,\tilde W}^*$ is
nonspecial.  So we expect (but have not proven!) that
for general $\tilde E$ a pullback of $E$ on $C$,
$K_{\tilde E,\tilde W}^*$ is nonspecial.
Proving this would prove Conjecture 2 and hence Conjecture 1.
Because then $K_{E,W}\otimes L_i$ would be nonspecial for
all $i$.

To prove this last conjecture it is required to prove that
a general $\tilde E$ is spanned.  That much we can do.
But unfortunately we work over $C$ and not $\tilde C$.
A new proof of the Second Theorem may lead to a proof
of the above.

Now we prove the Second Theorem.
\demo{Proof of Second Theorem}
A cyclic cover of degree $e$ (where $p$ does not divide $e$)
corresponds to a subgroup of $e$ torsion points
of the Jacobian.  Or in other words a set:
$$\{L, L^{\otimes 2} \dots L^{\otimes(e-1)}, L^{\otimes e} \},$$
with $L^{\otimes i} = \Cal O_C$ iff $i=e$.
Setting the above set of line bundles equal to 
${L_1, L_2, \dots L_e}$, and noting that the sections
of $\pi^* E$  are the pullback of $W$,
we then get the Second Theorem from Proposition 3.
\qed
\enddemo

\subheading{Remark 6}  Assuming the conjecture on the stability
of $K_{E,W}$ there are more equations (as in Remark 4) that
$g$ must solve.  We have $r=13$ and $d=5$ or $8$.  The case $e=2$ gives
nothing new (as in Remark 5).  But for $e=3$ we have 
$3\times 8 =24 = 13 + 11$ which gives the equation
$13(g-1)+8 \equiv 0 \pmod {11}$ or
$$g\equiv 8 \pmod {11}.$$  
But $e=3$ and $d = 5$ gives nothing new.
However, for $e=4$ and $d=5$ we get $20 =13 +7$ so 
$13(g-1)+5 \equiv 0 \pmod 7$.  
This reduces to 
$$g\equiv 6 \pmod 7.$$
Taking these equations, the previous equations, and using
the Chinese Remainder Theorem we get 
$$g \equiv 2,120 \pmod {2,310}.$$



\Refs

\ref \no 1 \by Ballico  \pages 21 -- 26
\paper Stable rationality for the  variety of vector
bundles over an algebraic curve
\jour J. Lond. Math. Soc. (2) \yr 1984 \vol 30
\endref

\ref \no 2  \by   A. Beauville, J-L. Colliot-Th\`el\'ene,
J-J Sansuc and P. Swinnerton-Dyer
\jour Ann. of Mathematics \yr 1985 \vol 121
\pages 283 -- 318
\paper Vari\'et\'es stablement  rationnelles
non rationnelles
\endref

\ref \no 3 
\by H.U. Boden and  K. Yokogawa 
\paper Rationality of moduli spaces of parabolic bundles
\jour Jour. of Lond. Math. Soc. \toappear
\year1995
\endref

\ref \no 4 \by H. Lange
\jour Jour. of Alg.
\paper Universal families of extensions
\year 1983 \vol 83 \pages 101 -- 112
\endref

\ref \no 5 \by M.S. Narasimhan and C.S. Seshadri
\paper Stable and unitary vector bundles on a 
compact Riemann surface
\jour  Ann. of Math. \vol 82 \yr 1969
\pages 540 -- 567
\endref


\ref \no 6 \by P.E. Newstead
\paper Rationality of moduli spaces of stable bundles
\yr 1975  \pages  251 -- 268 \vol 215
\jour Math. Ann.
\moreref \paper Correction
\yr 1980 \vol 249  \jour Math. Ann.  \pages 281 -- 282
\endref

\ref  \no 7 \by S. Ramanan 
\paper The moduli spaces of vector bundles  over an algebraic  curve
\jour Math. Ann. \vol 200 \yr 1973 \pages 69 -- 84
\endref

\ref \no 8 \by N.  Sundaram
\paper Special divisors and vector bundles
\yr 1987 \pages 175 -- 223 \vol  39
\jour Tohuko Math. J.
\endref

\ref \no 9 \manyby  A.N. Tyurin
\paper On the  classification of two dimensional
vector bundles over an  algebraic curve of
arbitrary genus
\yr 1964 \pages 21 -- 52 \vol 28
\jour Izv. Akad. Nauk. SSSR Ser.  Mat.
\endref

\ref \no 10 \bysame 
\paper Classification of vector bundles 
over an  algebraic curve of arbitrary genus
\yr 1964 \pages 657 -- 688 \vol 29
\jour Izv. Akad. Nauk. SSSR Ser.  Mat.
\moreref
\paper English translation
\jour  Amer. Math. Soc. Trans.
\vol 63 \yr 1967 \pages 245 -- 279
\endref


\endRefs
\enddocument